\documentclass[preprint,11pt,authoryear]{elsarticle}
\usepackage{graphicx}
\usepackage[linesnumbered,boxed]{algorithm2e}
\usepackage{mathrsfs}
\usepackage{amssymb}
\usepackage{subfigure}
\usepackage{hyperref}
\usepackage{multirow}
\usepackage{rotating}
\usepackage{hyperref}

\input amssym.def

\begin{document}

\renewcommand{\baselinestretch}{1.5}\large\normalsize

\begin{frontmatter}
\title{Iterated Tabu Search Algorithm for Packing Unequal Circles in a Circle}
\author[a,b]{Tao Ye}
\author[a]{Wenqi Huang}
\author[a]{Zhipeng L\"{u}\corref{cor1}}
\cortext[cor1]{Corresponding Author. Tel: 86-27-87543885. Email address: yeetao@gmail.com (T. Ye), wqhuang@hust.edu.cn(W. Huang), zhipeng.lui@gmail.com (Z. L\"{u})}
\address[a]{School of Computer Science and Technology, Huazhong University of Science and Technology, Wuhan, 430074, China}
\address[b]{Department of Mathematics, Simon Fraser University Surrey, Central City, 250-13450 102nd AV, Surrey, British Columbia, V3T 0A3, Canada }
\begin{abstract}
This paper presents an Iterated Tabu Search algorithm  (denoted by ITS-PUCC) for solving the problem of Packing Unequal Circles in a Circle.  The algorithm exploits the continuous and combinatorial nature of the unequal circles packing problem. It uses a continuous local optimization method to generate locally optimal packings. Meanwhile, it builds a neighborhood structure on the set of local minimum via two appropriate perturbation moves and integrates two combinatorial optimization methods,  Tabu Search and Iterated Local Search, to systematically search for good local minima. Computational experiments on two sets of widely-used test instances prove its effectiveness and efficiency. For the first set of  46  instances coming from the famous circle packing contest and the second set of 24 instances widely used in the literature, the algorithm is able to  discover respectively 14 and 16 better solutions than the previous best-known records.
\end{abstract}

\end{frontmatter}

\section{Introduction}
Given $n$ circles and a container of predetermined shape, the circle packing problem is concerned with a dense packing solution, which can pack all the circles into the smallest container without overlap. The circle packing problem is a well-known challenge in discrete and computational geometry, and  it arises in various real-world applications in the field of packing, cutting, container loading, communication networks and facility layout \citep{castillo2008}.  In the field of global optimization, the circle packing problem is a natural and challenging test bed for evaluating various global optimization methods.

This paper focuses on  solving  a classic  circle packing problem, the Packing Unequal Circles in a Circle (PUCC) problem.   As indicated in previous papers \citep{addis2008u, grosso2010, hifi2009},  the PUCC problem has an interesting and important characteristic that it has a both continuous and combinatorial nature. It has continuous nature because the position of each circle is chosen in $R^2$. The combinatorial nature is due to  the following two facts: (1) A packing pattern is composed of $n$ circles, and shifting a circle to a different place would produce a new packing pattern;  (2) The circles have different radiuses, and  swapping the positions of two different circles  may result in a new packing pattern.

 In this paper, we pay special attention to the continuous and combinatorial characteristic  of the PUCC problem. We propose an algorithm which integrates two kinds of optimization techniques: A continuous local optimization procedure which minimizes overlaps between circles and produces locally optimal packing patterns, and an Iterated Tabu Search (ITS) procedure which exploits the combinatorial nature of the problem and intelligently uses  two appropriate perturbation moves to search for globally optimal packing patterns.

The proposed algorithm is assessed on two sets of widely used test instances, showing its effectiveness and efficiency. For the first set of 46 instances coming from the famous circle packing contest, the algorithm is able to discover 14 better solutions than the previous best-known records. For the second set of 24 instances widely used in the literature, the algorithm can improve 16 best-known solutions in a reasonable time.

The rest of this paper is organized as follows. Section 2 briefly reviews the most related literature.  Section 3 formulates the PUCC problem. Section 4 describes the details of  the proposed algorithm. Section 5 assesses the performance of the  algorithm through extensive computational experiments.   Section 6 analyzes some key ingredients of the algorithm to understand the source of its performance. Finally, Section 7 concludes this paper and proposes some suggestions for future work.

\section{Related Literature}
Over the last few decades, the circle packing problem has received considerable attention in the literature. The simplest and most widely studied cases are the packing of equal circles in a square or in a circle.  Though researchers have spent significant effort on the two problems, only a few packings (up to tens of circles) have been proved to be optimal by purely analytical methods and computer-aided proving methods \citep{szabo2007,graham1998}. A second category of research aims at finding the best possible packings without optimality proofs. Following this spirit, various heuristic approaches have been proposed, including: Billiard simulation \citep{graham1998}, minimization of energy function \citep{nurmela1997}, nonlinear programming  approaches \citep{birgin2005, birgin2010}, Population Basin Hopping method \citep{addis2008e, grosso2010}, formulation space search heuristic algorithm \citep{beasley2011}, quasi-physical global optimization method \citep{huang2011}, greedy vacancy search method \citep{huang2010} and so on. With these approaches, best-known packings for up to thousands of circles have been found, which are reported  and continuously updated on the Packomania website \citep{specht2012}.

There are also a number of papers devoted to the unequal circle packing problem. Most previous papers on the unequal circle packing problem can be  classified into two categories: Constructive approaches and global optimization approaches. The constructive approaches build a packing by successively placing a circle into the container. These approaches usually include two important components: A placement heuristic, which determines several candidate positions for a new circle in the container, and a tree search strategy, which controls the tree search process and avoids exhaustive enumeration of  the solution space. The widely used placement heuristics include the principle of Best Local Position (BLP) \citep{hifi2004, hifi2007, hifi2008, akeb2010a} and the Maximum Hole Degree (MHD) rule \citep{huang2005, huang2006, lv2008, akeb2009}.  The tree search strategies include the self look-ahead search strategy \citep{huang2005, huang2006}, Pruned-Enriched-Rosenbluth Method (PERM) \citep{lv2008},   beam search algorithm \citep{akeb2009} and the hybrid beam search looking-ahead algorithm \citep{akeb2010a}.

The global optimization approaches formulate the unequal circle packing problem as a mathematical programming problem, then the task becomes to find  the global minimum of a mathematical model. These kind of approaches include the quasi-physical quasi-human algorithm by \citet{wang2002}, the Tabu Search and Simulated Annealing hybrid approach \citep{zhang2005}, the Population Basin Hopping algorithm  \citep{addis2008u, grosso2007, grosso2010}, the GP-TS algorithm by \citet{huang2012a}, the Iterated Local Search algorithm by \citet{huang2012b}, the Formulation Search Space algorithm by \citet{beasley2012} and the Iterated Tabu Search algorithm by \cite{fu2013} for the circular open dimension problem.

For the circle packing problem, there also exist many important literature not mentioned here. Interested readers are referred to the review articles by \citet{castillo2008} and  \citet{hifi2009},  the book by \citet{szabo2007} and the Packomania website \citep{specht2012}.

\section{Problem Formulation}
Given $n$ disks, each having radius $r_i$ $(i=1, 2, \cdots, n)$, the PUCC problem consists in finding a dense packing solution, which can pack all $n$ disks into the smallest circular container of radius $R$ without overlap. We designate the container center  as the origin of the cartesian coordinate system and locate disk $i$ $(i=1,2,\dots,n)$ by the coordinate position of its center $(x_i, y_i)$. The PUCC problem can be formulated as:
\begin{eqnarray}
 minimize \quad R ,  \quad s.t. : \nonumber   \\
    \sqrt{x_i^2 + y_i^2} + r_i \leq R  \\
    \sqrt{(x_i - x_j)^2 + (y_i - y_j)^2} \geq  r_i + r_j
\end{eqnarray}
where $ i,j = 1,2,\cdots,n; i \neq j$. Eq.(1) ensures that each disk is completely in the container and Eq.(2) guarantees that no overlap exists between any two disks. Note that, this problem can also be formulated in other ways, see for example \citet{birgin2005} and \citet{grosso2010}.

A packing solution is described by two variables: The radius of the container $R$ and the \textbf{packing pattern}  denoted by the positions of all $n$ disks $ X= (x_1, y_1,x_2, y_2, \cdots, x_n, y_n) $. The infeasibility of a packing  can be caused by two kinds of overlaps: Overlaps between two disks and overlaps between a disk and the exterior of the container. We define the overlapping depth between disks $i$ and $j$ $(i,j=1,2,\cdots,n;i\neq j )$ as:
\begin{equation}
\label{eq:oij}
	o_{ij}= max \ \{ 0 , r_i + r_j - \sqrt{(x_i - x_j)^2 + (y_i - y_j)^2} \} .
\end{equation}
and the overlapping depth between disk $i$ $( i=1,2,\cdots,n )$ and the exterior of the container as:
\begin{equation}
\label{eq:oxi}
	o_{0i}= max \ \{ 0, \sqrt{x_i^2 + y_i^2}+ r_i - R \} .
\end{equation}

Adding all squares of overlapping depth together, we get a penalty function  measuring overlaps of a packing
\begin{equation}
\label{eq:exs}
	E(X, R) = \sum_{i=0}^{n-1}{\sum_{j=i+1}^{n} o_{ij}^2} .
\end{equation}
Thus, a packing $(X, R)$ is  feasible (non-overlapping) if and only if $E(X, R)=0$.

 Sometimes, we fix the radius of the container at a constant value $R$ and the penalty function becomes
\begin{equation}
\label{eq:exs0}
	E_{R}(X) =  E(X, R).
\end{equation}
Note that, finding a packing pattern $X$ with $E_R(X)=0$ corresponds to solving the following circle-packing decision problem \citep{birgin2005}: Given a circular container with fixed radius $R$, find out a feasible pattern $X$ which can pack all the circles into the container without overlap.

Our original PUCC problem aims to find the smallest container of radius $R^*$ and a corresponding non-overlapping packing pattern $X$. In practice, the PUCC problem can be solved as a serial of circle-packing decision problems with descending $R$ \citep{huang2010}. The main steps are as follows:
\begin{enumerate}[(1)]
\item Let $\overline{R}$ be an upper bound of $R^*$. Initialize $\overline{R}$ with a relatively large number such that all circles can be easily packed into the container of radius $\overline{R}$ without overlap.
\item Set $R \leftarrow \overline{R}$ and launch an algorithm to find a feasible $X$ with $E_R(X) = 0$ (i.e., to solve the corresponding circle-packing decision problem).

\item Tighten the packing $(X, R)$, i.e, to minimize $R$ while keeping $X$ basically unchanged. This step can be achieved using various approaches, like the simple bisection method described in \cite{huang2011},  the simple penalty method described in \cite{huang2010},  the standard local optimization solver SNOPT adopted by \cite{addis2008u} and the more sophisticated augmented Lagrangian method \citep{andreani2007, birgin2008,birgin2009}. After this step, we can usually obtain a better (at least not worse) packing $(X', R')$.
\item  Set $\overline{R} \leftarrow R'$ and go to step 2. The loop of steps 2-4 is ended until a certain termination criterion (like time limit) is satisfied.
\end{enumerate}

In the rest of this paper, we will first introduce an Iterated Tabu Search algorithm to solve the circle-packing decision problem, and then use it to search for dense packing solutions for the PUCC problem in the computational experiments section.

\section{Iterated Tabu Search Algorithm}
\label{sec:ITS}
This section describes the Iterated Tabu Search (ITS) procedure for solving the circle-packing decision problem. As indicated in Section 3, this problem can be transformed to an unconstraint global optimization problem:
\begin{equation}
	 minimize \ \ \ E_R(X).
\end{equation}
This subproblem is very difficult because there exist enormous local minima in the solution space. \citet{grosso2010} have shown that, even for the equal circle packing problem, the number of local minima tends to increase very quickly with the number of circles. For the more complex unequal circle packing problem, it is very possible that the number of local minima will be significantly larger.

The main rationale behind the ITS procedure is as follows: (1) Each local minimum of $E_R(X)$ corresponds to a packing pattern of $n$ disks in the container. (2) If we perturb the current local minimizer $X$ by swapping the positions of two different disks (or shifting the position of one disk) and then call the LBFGS procedure to minimize $E_R(X)$, we can obtain a new local minimizer. (3) By systematically using the two perturbation moves, swap and shift, we can obtain a set of neighboring local minima from the current local minima.  Furthermore, we can build a neighborhood structure on the set of local minima of $E_R(X)$.  (4) Since there is a neighborhood structure,  some  Stochastic Local Search methods \citep{slsbook}, such as Tabu Search \citep{tabusearch} and Iterated Local Search \citep{ils},  can be employed to search for good local minima.

The outline of the ITS procedure is given in Algorithm \ref{alg:its}.  The procedure performs searches on the set of  local minima of $E_R(X)$ and follows an Iterated Local Search schema.  In Algorithm \ref{alg:its}, we run the ITS procedure  in a multi-start fashion. At each run, the algorithm starts from a randomly generated local minimum (steps 2-3). It goes through the $SwapTabuSearch$ procedure (step 4) and reaches a swap-optimal local minimum (which will be defined in the next section).  Then the search explores the solution space by repeatedly escaping from  local optima traps (step 6) and moving to another local optimum (step 7). This process is repeated until the best-found solution has not been improved during the last $PerturbDepth$ iterations.

\begin{algorithm}
\label{alg:its}
\caption{The ITS procedure}
\KwIn{Radiuses of $n$ disks, radius of the container $R$}
\KwOut{A feasible packing pattern $X$ of $n$ disks in the container }
 \While {not time out}
 {
    $X \leftarrow $ Randomly scatter $n$ disks into the container \;
    $X \leftarrow $ Minimize $E_R(X)$ using LBFGS \;
    $X \leftarrow $ SwapTabuSearch($X$)  \tcc*{local search}
    \Repeat{the best-found solution has not been improved in the last PerturbDepth iterations }
    {
        $Y \leftarrow$ ShiftPerturb ($X$) \tcc*{perturb}
        $Y \leftarrow $ SwapTabuSearch($Y$)  \tcc*{local search}
        \If{$E_R(Y) \leq E_R(X)$}
        {
            $X \leftarrow Y$;
        }
    }
 }
\end{algorithm}

\subsection{The SwapTabuSearch Procedure}
\label{sec:tabusearch}
In the $SwapTabuSearch$ procedure, we build a neighborhood structure on the set of local minima of $E_R(X)$.   A \textbf{swap move} performed on a packing pattern $X$ is defined as swapping the positions of two disks with different radiuses and then locally minimizing $E_R(X)$. For two local minima of $E_R(X)$, $X$ and $X'$, we say $X'$  is a \textbf{neighbor} of $X$ if and only if $X'$  can be reached by performing a swap move on $X$. The \textbf{neighborhood} of $X$ is  a set containing all the neighbors of $X$. We use $E_R(X)$ as the evaluation function, and a local minimum $X$ of $E_R(X)$ is called a \textbf{swap-optimal} local minimum if it has better solution quality than all its neighbors (or it cannot be improved via any swap move).

Totally, there are $n*(n-1)$ possible swap moves for a packing pattern with $n$ disks. However, for efficiency purposes, a restricted neighborhood is used in this paper. We first sort the disks in a nondecreasing order w.r.t. their radius values, such that for disks $i$ and $j $,   $r_i \leq r_j$ if $i < j$.  A swap move can only be performed on a pair of disks with neighboring radius values. That is to say, disk $i$ can only exchange positions with disks $i-1$ and $i+1$. Then, there are in total $n-1$ swap moves and a local minimum $X$ has at most  $n-1$ different neighbors.

The $SwapTabuSearch$ procedure follows a Tabu Search strategy. At the beginning of the search, the tabu list is empty and all swap moves are admissible. At each step,  the algorithm chooses a best admissible move which leads to the best nontabu solution. The aspiration criterion is used such that a tabu move can be selected if it generates a solution that is better than the best-found solution.  Once a move is selected, it is declared tabu for the next $TabuTenure$ steps. The procedure is repeated until the best-found solution has not been improved with the last $TabuDepth$ steps. The sketch of the $SwapTabuSearch$ procedure is presented in Algorithm \ref{alg:SwapTabuSearch}.

\begin{algorithm}
\label{alg:SwapTabuSearch}
\caption{The SwapTabuSearch procedure}
\KwIn{An initial local minimum of $E_R(X)$}
\KwOut{A swap-optimal minimum of $E_R(X)$ }
 \Repeat {The best-found solution has not been improved in the last $TabuDepth$ iterations}
 {
       choose a best candidate swap move $mv$\;
       perform move $mv$ \;
       declare move $mv$ tabu for $TabuTenure$ iterations\;
 }
 \Return the best-found solution $X^*$ \;
\end{algorithm}

\subsection{The ShiftPerturb Procedure}
\label{sec:perturb}
In the $ShiftPerturb$ procedure, the algorithm escapes a local optimum by a series of shift moves. A shift move performed on a packing pattern $X$ is defined as shifting the position of  a randomly chosen disk to a random place in the container and then locally minimizing $E_R(X)$. The number of times the shift move is performed is controlled by a parameter $PerturbStrength$. As pointed out in previous research \citep{ils}, the perturbation strength is very important for Iterated Local Search. If it is too weak, the local search may undo the perturbation and the search will be confined in a small area of the solution space. On the contrary, if the perturbation is too strong, the Iterated Local Search will behave like random restart, leading to poor performance. After preliminary computational tests, we choose the value of $PerturbStrength$ to be a random integer from $[1, n/8]$.

\section{Performance Assessment}
In this section, we assess the performance of the proposed algorithm through computational experiments on two sets of widely-used test instances. We also compare the results of our algorithm with some state-of-the-art algorithms in the literature.

\subsection{Experimental Protocol}
The algorithm is programmed in C++ and complied using GNU G++. All  computational experiments are carried out on a personal computer with 4Gb memory and a 2.8GHz AMD Phenom II X6 1055T CPU. Table \ref{tbl:parameter} gives the settings of the four  important parameters of the algorithm. Note that all the computational results are obtained without special tuning of the parameters, i.e., all the parameters used in the algorithm are fixed for all the tested instances.

\begin{table}
\begin{center}
\caption{Parameter settings}
\label{tbl:parameter}
\scriptsize{
\begin{tabular}{|p{2.2cm}p{1cm}p{4.5cm}p{2.3cm}|}
\hline
Parameter               & Section  &     Description & Value  \\
\hline
$TabuTenure$          &  \ref{sec:tabusearch} &    Tabu tenure of Tabu Search  &$n/5$ + rand(0,10)  \\
$TabuDepth $          & \ref{sec:tabusearch} &    Improvement cutoff of Tabu Search &10*$n$         \\
$PerturbStrength$    & \ref{sec:perturb} &   Perturbation strength of ITS & rand(1,$n$/8)         \\
$PerturbDepth$        & \ref{sec:ITS} &    Improvement cutoff of ITS  &10*$n$         \\
\hline
\end{tabular}
}
\end{center}
\end{table}

\subsection{Test Instances}
\label{sec:instances}
Two sets of test problems are considered, in total constituting 70 instances. The first set comes from the famous circle packing contest (see  \url{http://www.recmath.org/contest/CirclePacking/index.php}). This contest started on October 2005 and ended on January 2006. During this period, the participants were invited to propose densest packing solutions to pack $n(n=5,6,\dots,50)$ circles, each having radius $r_i = i (i=1, 2,\dots, n)$ into the smallest containing circle without overlap.  155 groups from 32 countries took part in the contest and  submitted a total of 27490 tentative solutions. After the contest, these results were further improved respectively by \citet{muller2009}, Eckard Specht \citep{specht2012}, Zhanghua Fu et al. \citep{specht2012}. Currently, all the best-known records are published and continuously updated on the Packomania website.

The second set of instances consists of 24 problem instances first presented by \cite{huang2005}. These instances are frequently used in the literature by many authors, see for example \cite{huang2006, akeb2009, akeb2010a}. The size of these instances ranges from $n=10$ to $60$. A detailed description of these instances can be found in \cite{huang2005}.

\subsection{Computational results on the circle packing contest instances}
For the circle packing contest instances,  researchers usually pay more attention on the solution quality. Especially during the contest,  people mostly focus on finding better solutions than the best-known records and rarely consider the computational resource used. After researchers have solved these instances using various approaches and large amount of computational resource, this set of instances now becomes a challenging benchmark to test the \textit{discovery capability} \citep{grosso2007} of a new algorithm.  Therefore, our first  experiment concentrates on searching for high-quality solutions. For each run of each instance, we usually set the time limit to 24 hours, run the algorithm multiple times and record the best-found solutions.

Table \ref{tbl:ccin_result} gives the computational results. Column 1 lists the best-known records on the Packomania website. Columns 2-4 respectively report the solution difference between some top reference results and the best-known records. These include: the best results found by \cite{addis2008u} (who is the champion of the circle packing contest) using PBH algorithm, the best records obtained by all the participants in the contest, the best results found by \cite{muller2009} using Simulated Annealing (SA) algorithm. Column 5 gives the solution difference between our results and  the  best-known records. The results indicated in bold are better than the best-known ones. Table \ref{tbl:ccin_result}  omits  the results for  $n=5, 6, \dots, 20$, because our results and all the reference results are the same on these instances.  Note that, our program generates solutions with a maximum error on the distances of $10^{-9}$. We have sent all the improved results to Eckard Specht. Using his own local optimization solver, he has processed our results to a high precision ($10^{-28}$) and published them on the Packomania website.

Table \ref{tbl:compare} summarizes the comparison of our results with the reference results. The rows \textit{better}, \textit{equal} and \textit{worse} respectively denote the number of instances for which the proposed algorithm gets solutions that are better, equal and worse than each reference result.  Table \ref{tbl:compare} shows that the proposed algorithm is able to discover a number of better solutions than the previous best reference results, demonstrating its efficacy in finding high-quality solutions. In fact, we also tested the proposed algorithm on the larger instances of $n=51, 52, \dots, 100$. Some preliminary experiments show that the algorithm can improve almost all previous best-known results. Interested readers can refer to the Packomania website.

\renewcommand{\baselinestretch}{1.0}\large\normalsize
\begin{table}
\begin{center}
\caption{Comparison of solution quality on the circle packing contest instances}
\label{tbl:ccin_result}
\scriptsize{
\begin{tabular}{cccccc}
\hline
\multirow{2}{*}{$n$} 	&	\multirow{2}{*}{Best-Known}	 &\multicolumn{4}{c}{Solution difference (i.e., this result  - best-known) } \\
\cline{3-6}
   &   &	 PBH	& Contest record	&	 SA	&	ITS-PUCC	\\
\hline
21	&	62.55887709 	&	0.00118149 	&	0	&	0	&	0	\\
22	&	66.76028624 	&	0	&	0	&	0	&	0	\\
23	&	71.19946160 	&	0	&	0	&	0	&	0	\\
24	&	75.74914258 	&	0.00356154 	&	0.00356154 	&	0	&	0	\\
25	&	80.28586443 	&	0	&	0	&	0	&	0	\\
26	&	84.98993916 	&	0.11634365 	&	0.08646206 	&	0	&	\textbf{-0.01174810} 	\\
27	&	89.75096268 	&	0.07861113 	&	0.04121888 	&	0	&	0	\\
28	&	94.52587710 	&	0.17998508 	&	0.02410937 	&	0.0006594	&	0	\\
29	&	99.48311156 	&	0.02920634 	&	0.02920634 	&	0	&	0	\\
30	&	104.54036376 	&	0.20743552 	&	0.03819132 	&	0.0008052	&	0	\\
31	&	109.68204275 	&	0.08990423 	&	0.08990423 	&	0.0004137	&	\textbf{-0.05280209} 	\\
32	&	114.79981466 	&	0.06562367 	&	0.06562367 	&	0.0411343	&	0	\\
33	&	120.06565963 	&	0.15129751 	&	0.15129751 	&	0.0001869	&	0	\\
34	&	125.36693920 	&	0.24548366 	&	0.06656255 	&	0	&	0.07661871 	\\
35	&	130.84907874 	&	0.31742394 	&	0.30727589 	&	0.0685492	&	0	\\
36	&	136.49212355 	&	0.04277728 	&	0.04277728 	&	0.0001210	&	\textbf{-0.18421273} 	\\
37	&	141.93243775 	&	0.32389631 	&	0.24254278 	&	0.1189377	&	\textbf{-0.14870433}	\\
38	&	147.45211646 	&	0.53945928 	&	0.40557489 	&	0.0047288	&	0.16654317 	\\
39	&	153.30070280 	&	0.30312533 	&	0.25459839 	&	0.0793799	&	\textbf{-0.00009525} 	\\
40	&	159.17977260 	&	0.39413040 	&	0.30925227 	&	0.0026352	&	\textbf{-0.12653569} 	\\
41	&	164.88704217 	&	0.40486751 	&	0.40486751 	&	0.1498584	&	0	\\
42	&	170.89531908 	&	0.03044253 	&	0.03044253 	&	0.0000083	&	\textbf{-0.11479840} 	\\
43	&	176.82574386 	&	0.41388066 	&	0.24859621 	&	0.2256308	&	\textbf{-0.02750344} 	\\
44	&	183.04328935 	&	0.32678190 	&	0.13277222 	&	0.0559354	&	\textbf{-0.11226929} 	\\
45	&	189.19513856 	&	0.48403531 	&	0.44030054 	&	0.0077934	&	\textbf{-0.02405842} 	\\
46	&	195.52636407 	&	0.38439932 	&	0.38439932 	&	0.0000710	&	\textbf{-0.17039157} 	\\
47	&	201.72792559 	&	0.50009381 	&	0.45768615 	&	0	&	\textbf{-0.05939136} 	\\
48	&	208.09015930 	&	0.54578742 	&	0.54578742 	&	0	&	\textbf{-0.04513139} 	\\
49	&	214.29205550 	&	0.36989651 	&	0.36989651 	&	0.0033920	&	\textbf{-0.00983298} 	\\
50	&	220.56540026 	&	0.52435233 	&	0.52435233 	&	0.0350184	&	0.50085590 	\\
\hline
\end{tabular}
}
\end{center}
\end{table}
\renewcommand{\baselinestretch}{1.5}\large\normalsize

\renewcommand{\baselinestretch}{1.0}\large\normalsize
\begin{table}
\begin{center}
\caption{Summary of comparison of solution quality on the 30 circle packing contest instances with $n=21, 22, \dots, 50$}
\label{tbl:compare}
\scriptsize{
\begin{tabular}{|c|cccc|}
\hline
              & PBH       & Contest record       & SA  & Best-known \\
\hline
Better    &   27                   &    25                        & 20                         &  14                    \\
Equal    &    3                    &     4                         &   7                         &  13                   \\
Worse   &    0                    &     1                         &   3                         &   3                   \\
\hline
\end{tabular}
}
\end{center}
\end{table}
\renewcommand{\baselinestretch}{1.5}\large\normalsize

All the reference algorithms in Table \ref{tbl:ccin_result} concentrate on finding high-quality solutions and do not reveal their computational statistics.  In order to further evaluate the proposed algorithm in terms of search efficiency, we conduct additional experiments to compare the proposed algorithm with two recently published algorithms in a time-equalized basis. For each instance of $n=5,6, \dots, 32$, we set the maximum time limit to 10000 seconds.  We record the best-found solution and the elapsed time when it is first detected by the algorithm. To reduce the impact of randomness, each instance is independently solved for 10 times.

Table \ref{tbl:contest2} gives the computational results. Columns 2-3, 4-5 respectively list the best-found solution and the needed computing time of TS/NP algorithm and  FSS algorithm.  Columns 2 and 3 are extracted from \cite{aimudahka2010} where the algorithm ran on a computer with a Pentium IV, 2.66 Ghz CPU and 512Mb RAM. Columns 4 and 5 are extracted from \cite{beasley2012}. Their experiments were done on a computer with a Intel(R) Core(TM) i5-2500 3.30 GHz CPU and 4.00 GB RAM. Columns 6-8 give the computational statistics of our algorithm, including the best-found solution, the number of hit times and the averaged computing time to detect the best-found solution.

Columns 6-8 show that, for all the 28 instances, the proposed algorithm can reach (or improve) the previous best-known records listed in Table \ref{tbl:ccin_result} within the given time limit. Especially for $n\leq 25$, the algorithm can robustly detect the best-known records in a short time. When compared with the two reference algorithms, one observes that  the proposed algorithm can usually find better solutions within the time limit. These results provide evidence of  the search efficiency of  ITS-PUCC algorithm.
\renewcommand{\baselinestretch}{1.0}\large\normalsize
\begin{table}
\begin{center}
\caption{Comparison of search efficiency on the circle packing contest instances}
\label{tbl:contest2}
\scriptsize{
\begin{tabular}{p{0.5cm}p{1.6cm}p{0.8cm}p{0.01cm}p{1.6cm}p{0.8cm}p{0.01cm}p{1.6cm}p{0.8cm}p{0.8cm}}
\hline
\multirow{2}{*}{$n$} & \multicolumn{2}{c}{TS/NP}& & \multicolumn{2}{c}{FSS} & &\multicolumn{3}{c}{ITS-PUCC} \\
\cline{2-3}  \cline{5-6}  \cline{8-10}
  & $R$ & time(s) & & $R$ & time(s) & & $R$ & \#hits & time(s)  \\
\hline
5	&	9.001398	&	426	& & 	9.00139775	&	461 	& &	9.00139774	&	 10/10	&	1	\\
6	&	11.05704	&	686	& & 	11.0570404	&	667 	& &	11.05704039	&	 10/10	&	1	\\
7	&	13.46211	&	1511	& & 	13.46211068	&	721 	& &	13.46211067	&	 10/10	&	1	\\
8	&	16.22175	&	2551	& & 	16.22174668	&	1028 	& &	16.22174667	&	 10/10	&	1	\\
9	&	19.39734	&	4051	& & 	19.23319391	&	1404 	& &	19.2331939	&	 10/10	&	1	\\
10	&	22.34516	&	5760	& & 	22.00019301	&	1438 	& &	22.00019301	&	 10/10	&	1	\\
11	&	24.96063	&	2094	& & 	24.96063429	&	1820 	& &	24.96063428	&	 10/10	&	1	\\
12	&	28.67863	&	3548	& & 	28.37138944	&	2299 	& &	28.37138943	&	 10/10	&	1	\\
13	&	32.00719	&	4054	& & 	31.54586702	&	2905 	& &	31.54586701	&	 10/10	&	1	\\
14	&	35.41261	&	6146	& & 	35.09564714	&	2970 	& &	35.09564714	&	 10/10	&	2	\\
15	&	39.00243	&	6696	& & 	38.83799682	&	3904 	& &	38.8379955	&	 10/10	&	1	\\
16	&	42.92185	&	9684	& & 	42.45811644	&	4917 	& &	42.45811643	&	 10/10	&	6	\\
17	&	46.77237	&	10168	& & 	46.34518193	&	5264 	& &	46.29134211	&	 10/10	&	24	\\
18	&	50.65635	&	14312	& & 	50.20889346	&	6224 	& &	50.11976262	&	 10/10	&	23	\\
19	&	55.02744	&	14925	& & 	54.36009421	&	7349 	& &	54.24029359	&	 10/10	&	39	\\
20	&	59.04547	&	19825	& & 	58.48047359	&	7517 	& &	58.40056747	&	 10/10	&	82	\\
21	&	63.49768	&	5923	& & 	63.00078332	&	8924 	& &	62.55887709	&	 10/10	&	190	\\
22	&	68.10291	&	6636	& & 	66.96471591	&	10762 	& &	66.76028624	&	 10/10	&	127	\\
23	&	72.70501	&	7209	& & 	71.69822657	&	13018 	& &	71.1994616	&	 10/10	&	268	\\
24	&	76.49105	&	8552	& & 	76.1231197	&	13004 	& &	75.74914258	&	 10/10	&	704	\\
25	&	81.56595	&	11409	& & 	80.8168236	&	15569 	& &	80.28586443	&	 10/10	&	633	\\
26	&	86.43809	&	12062	& & 	85.487438	&	18320 	& &	84.97819106	&	 9/10	&	3538	\\
27	&	91.15366	&	13657	& & 	90.93173506	&	18544 	& &	89.75096268	&	 7/10	&	5287	\\
28	&	96.34813	&	14364	& & 	95.6406414	&	21931 	& &	94.5258771	&	 10/10	&	1568	\\
29	&	101.7251	&	15185	& & 	100.7200313	&	25455 	& &	99.48311156	&	 7/10	&	2915	\\
30	&	107.1161	&	20745	& & 	105.8881722	&	25658 	& &	104.5403638	&	 5/10	&	4538	\\
31	&	111.8996	&	21424	& & 	111.077126	&	29973 	& &	109.6292407	&	 2/10	&	8551	\\
32	&	117.6701	&	22781	& & 	116.6122668	&	34445 	& &	114.7998147	&	 4/10	&	3885	\\
\hline
\end{tabular}
}
\end{center}
\end{table}
\renewcommand{\baselinestretch}{1.5}\large\normalsize
\subsection{Computational results on the  NR instances}
This section tests the proposed algorithm on the 24 NR instances. For each instance, we set the time limit to 10000 seconds, and record the best-found solution and the elapsed time when it is first detected by the algorithm. Each instance is solved for 10 times from different randomly generated starting points.

The computational results are presented in Table \ref{tbl:NR}. Column 1 gives the instance name. Columns 2-3, 4-5, 6-7, 8-9, respectively present the best-found solution and the needed computing time of  A1.5 Algorithm in \cite{huang2006}, Beam Search (BS) algorithm in \cite{akeb2009}, Algorithm 2 in \cite{akeb2010b} and GP-TS algorithm in \cite{huang2012a}. Columns 10-12 give the computational statistics of our algorithm, including the best-found solution, the number of hit times and the averaged computing time for detecting the best-found solution. In experiments, our  program generates solutions with a maximum error on the distance of $10^{-9}$. However, in order to keep consistent with previous papers, we report in Table \ref{tbl:NR} the results with 4 significant digits after the decimal point.

Table \ref{tbl:NR} demonstrates that, for all the tested 24 instances, the proposed algorithm can find 16 better solutions than the best results found by the references algorithms (as indicated in bold in the table). For the other 8 instances, it can reach the best-known solutions efficiently and robustly. These results further provide evidence of the competitiveness of the proposed algorithm.
\renewcommand{\baselinestretch}{1.0}\large\normalsize
\begin{sidewaystable}
\begin{center}
\caption{Comparison of search efficiency on the 24 NR instances}
\label{tbl:NR}
\scriptsize{
\begin{tabular}{p{1.0cm}p{0.6cm}p{0.3cm}p{0.01cm}p{1.0cm}p{0.3cm}p{0.01cm}p{1.0cm}p{0.4cm}p{0.01cm}p{1.0cm}p{0.3cm}p{0.01cm}p{1.1cm}p{0.8cm}p{0.6cm}}
\hline
\multirow{2}{*}{Instance} & \multicolumn{2}{c}{A1.5}& & \multicolumn{2}{c}{BS} & &\multicolumn{2}{c}{Algorithm 2} && \multicolumn{2}{c}{GP-TS} && \multicolumn{3}{c}{ITS-PUCC}\\
\cline{2-3}  \cline{5-6}  \cline{8-9}  \cline{11-12}  \cline{14-16}
 & $R$ & time(s) & & $R$ & time(s) & & $R$ & time(s) && $R$ & time(s) && $R$ & \#hits & time(s) \\
\hline
NR10-1 & 99.89 & 1 && 99.8851 & 19 && 99.8851 & 1 && 99.8851 & 1 && 99.8851 & 10/10 & 1 \\
NR11-1 & 60.71 & 1 && 60.7100 & 28 && 60.7100 & 2 && 60.7100 & 2 && 60.7100 & 10/10 & 1 \\
NR12-1 & 65.30 & 6 && 65.4752 & 4 && 65.0338 & 2483 && 65.0245 & 1 && 65.0244 & 10/10 & 1 \\
NR14-1 & 113.84 & 2 && 114.2919 & 151 && 113.5588 & 17860 && 113.5588 & 252 && 113.5588 & 10/10 & 16 \\
NR15-1 & 38.97 & 25 && 38.9441 & 59 && 38.9170 & 46283 && 38.9158 & 88 && \textbf{38.9114} & 10/10 & 57 \\
NR15-2 & 38.85 & 6 && 38.8380 & 1179 && 38.8380 & 9832 && 38.8380 & 166 && 38.8380 & 10/10 & 1 \\
NR16-1 & 143.44 & 71 && 143.7176 & 139 && 143.4339 & 235240 && 143.3798 & 128 && 143.3798 & 10/10 & 126 \\
NR16-2 & 128.29 & 44 && 128.0539 & 28 && 127.9021 & 80890 && 127.7174 & 6 && \textbf{127.6978} & 10/10 & 2 \\
NR17-1 & 49.25 & 30 && 49.2069 & 234 && 49.1977 & 6080 && 49.1874 & 258 && 49.1873 & 10/10 & 359 \\
NR18-1 & 197.40 & 88 && 198.2791 & 8 && 197.1038 & 111970 && 197.0367 & 76 && \textbf{196.9826} & 10/10 & 58 \\
NR20-1 & 125.53 & 39 && 125.6316 & 764 && 125.1525 & 199945 && 125.1178 & 13 && 125.1178 & 10/10 & 229 \\
NR20-2 & 122.21 & 318 && 122.2192 & 351 && 122.0296 & 150495 && 121.9944 & 120 && \textbf{121.7887} & 10/10 & 780 \\
NR21-1 & 148.82 & 683 && 149.1351 & 638 && 148.3462 & 132080 && 148.3373 & 647 && \textbf{148.0968} & 10/10 & 488 \\
NR23-1 & 175.47 & 1229 && 175.4058 & 3072 && 174.9491 & 58430 && 174.8524 & 283 && \textbf{174.3425} & 8/10 & 3259 \\
NR24-1 & 138.38 & 2339 && 138.2778 & 510 && 138.0520 & 37140 && 138.0044 & 433 && \textbf{137.7591} & 10/10 & 3706 \\
NR25-1 & 190.47 & 4614 && 190.1855 & 1493 && 189.4715 & 37053 && 189.3736 & 533 && \textbf{188.8314} & 7/10 & 3091 \\
NR26-1 & 246.75 & 1019 && 247.5464 & 583 && 246.4179 & 18620 && 246.0853 & 40 && \textbf{244.5743} & 10/10 & 795 \\
NR26-2 & 303.38 & 5164 && 303.2102 & 11240 && 302.5896 & 190600 && 302.0687 & 4 && \textbf{300.2631} & 10/10 & 179 \\
NR27-1 & 222.58 & 4436 && 222.4896 & 3750 && 221.6389 & 177000 && 221.4882 & 2505 && \textbf{220.9393} & 5/10 & 6037 \\
NR30-1 & 178.66 & 1365 && 178.0102 & 5045 && 177.6473 & 160700 && 178.0093 & 10439 && \textbf{177.5125} & 3/10 & 8678 \\
NR30-2 & 173.70 & 1078 && 173.4359 & 9217 && 173.2215 & 155050 && 173.1641 & 502 && \textbf{172.9665} & 8/10 & 574 \\
NR40-1 & 357.00 & 12109 && 357.0695 & 22140 && 355.6587 & 158700 && 355.1307 & 112 && \textbf{352.4517} & 4/10 & 2253 \\
NR50-1 & 380.00 & 9717 && 378.5854 & 21400 && 378.0044 & 186000 && 377.9105 & 12462 &&\textbf{377.9080} & 2/10 & 8945 \\
NR60-1 & 522.93 & 13256 && 521.2739 & 13975 && 519.8494 & 116270 && 519.4515 & 60 &&\textbf{518.6792} & 1/10 & 9657 \\
\hline
\end{tabular}}
\end{center}
\end{sidewaystable}
\renewcommand{\baselinestretch}{1.5}\large\normalsize

\section{Algorithm Analysis}
In this section, we turn our attention to  analyzing the two most important ingredients of the proposed algorithm: the $SwapTabuSearch$ procedure and the $ShiftPerturb$ procedure.

\subsection{Analysis of The SwapTabuSearch Procedure}
The $SwapTabuSearch$ procedure is a key component of the proposed algorithm, which enables the algorithm to intelligently examines the neighboring packing patterns through swap moves. In order to make sure the Tabu Search strategy makes a meaningful contribution, we conduct experiments to compare the Tabu Search strategy with a simple local search strategy called Steepest Descent \citep{slsbook}.

For comparison, we use the same neighborhood structure as described in Section \ref{sec:tabusearch}  and implement the Steepest Descent strategy as follows. At each iteration, the search examines each neighbor of the current solution and find out the best neighbor with the least objective value $E_R$. If the best neighbor $X'$ is better than the current solution $X$, i.e., $E_R(X') \leq E_R(X)$, then the search moves to $X'$ and proceeds to the next iteration;  otherwise the search stops and declares reaching a local minimum.

A representative instance \textit{NR15-2} is chosen as our test bed. This instance is nontrivial. Though many previous papers have tested it, only  few state-of-the-art algorithms, like Beam Search \citep{hifi2008}, PBH\citep{addis2008u,grosso2010}, SA\citep{muller2009} can obtain the optimal packing pattern. We set the radius of container $R$ to the best-known value, randomly generate initial $X$ and call both algorithms to minimize $E_R(X)$.

We run both algorithms 1000 times from different randomly generated starting points and record in Table \ref{tbl:ts_sd} respectively the best-found solution (Column 2), the average solution quality (Column 3), the average number of search steps for each local search (Column 4) and the average elapsed time for each local search (Column 5).  From Table \ref{tbl:ts_sd}, we observe that, the Tabu Search strategy shows clear advantage over Steepest Descent strategy. Each time, the Tabu Search strategy can find the global minimum from a randomly generated starting point, while the Steepest Descent strategy fails for all 1000 runs. In fact, we try to run the Steepest Descent strategy from 100000 randomly generated starting points, it still cannot find the global minimum.

\begin{figure}
 \includegraphics[width=4in]{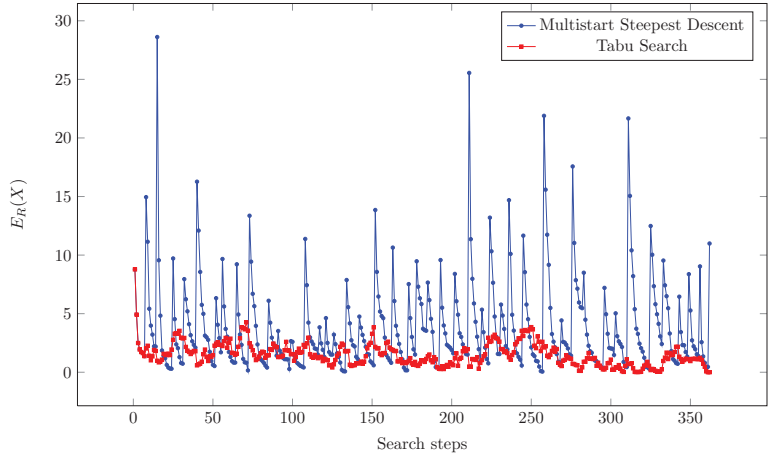}
 \caption{Comparison of search trajectories between  Multistart Steepest Descent and Tabu Search }
 \label{fig:sd_ts}
\end{figure}

The main reason for the difference is that, with  the Steepest Descent strategy, the search is easily trapped in poor local minimum. As shown in Table \ref{tbl:ts_sd}, the average number of search step for each local search is only 7. However, with the Tabu Search strategy, the search can escape from low-quality local minimum trap and proceed to explore the neighboring area. Figure \ref{fig:sd_ts} shows a typical search trajectory of Tabu Search, compared with the search trajectory of Multistart Steepest Descent. In Figure \ref{fig:sd_ts}, both algorithms start from the same initial solution, a packing pattern with $E_R = 8.77845794$. After 7 search steps, both of them encounter a local minimum with $E_R = 1.37297106 $. At this time, the Steepest Descent strategy is trapped,  the search has to proceed from a new randomly generated initial solution. However, with the Tabu Search strategy, the search is able to escape from the local minimum with $E_R = 1.3729106$, proceed to examine the neighboring area, and finally find the global minimum at the 362th search step.

\begin{table}
\begin{center}
\caption{Computational statistics of Tabu Search strategy and Steepest Descent strategy from 1000 randomly generated initial packings}
\label{tbl:ts_sd}
\scriptsize
\begin{tabular}{|c|cccc|}
\hline
  Search Strategy & Best-found solution & Average solution quality  &   Search steps  & Time (s) \\
 \hline
 Tabu Search & 0.000000 & 0.000000 & 1269 & 2 \\
 \hline
  Steepest Descent & 0.092194 & 2.085742 & 7 & 0 \\
  \hline
  \end{tabular}
  \end{center}
\end{table}

These experiments reveal that,  the Tabu Search strategy helps to perform an intensified examination around the incumbent packing pattern and  makes possible discovering those hidden good solutions.  The same experiments have been performed on several other instances, leading to similar observation.

\subsection{Analysis of the ShiftPerturb Procedure}
In order to verify the effectiveness of the $ShiftPerturb$ procedure, we conduct experiments to compare the  proposed ITS algorithm with a Multistart Tabu Search algorithm. In the Multistart Tabu Search algorithm, when the $SwapTabuSearch$ procedure finishes, the search proceeds from a new randomly generated initial solution. The parameter setting of $SwapTabuSearch$ is the same as listed in Table \ref{tbl:parameter}. We test the Multistart Tabu Search algorithm on the 28 circle packing instances with $5 \leq n \leq 32$. Each instance is solved for 10 times. The time limit for each run is also set to 10000 seconds.

The computational results show clear advantage of  ITS algorithm over Multistart Tabu Search algorithm. For the instances of $5 \leq n\leq 24$, the Multistart Tabu Search algorithm can also detect the best-known records, but with lower success rates and relatively longer time.  Nevertheless, for each instance of $25 \leq n \leq 32$, the Multistart Tabu Search algorithm fails to detect the best-known solution for all the 10 runs within the given time limit.

We conjecture the superior of ITS over Multistart Tabu Search may be explained from the following two aspects. First, the Iterated Local Search framework helps the search to perform a more intensified examination around the incumbent solution, making it possible to repeatedly discover better solutions. This is supported by our observations from computational experiments that, with the ITS algorithm, the search can usually generate a sequence of local minima with descending objective value. The final solution obtained by one run of ITS is usually much better than that found by the first run of $SwapTabuSearch$. Second, the shift move in the $ShiftPerturb$ procedure is complementary to the swap move, enabling the search to reach some packing patterns which are hard to detect only through swap moves.

\section{Conclusion and Future Work}
In this paper, we have presented a heuristic global optimization algorithm for solving the unequal circle packing problem. The proposed algorithm uses a continuous local optimization method to generate locally optimal packings and integrates two combinatorial optimization methods,  Tabu Search and Iterated Local Search, to systematically search for good local minima. The efficiency and effectiveness of the algorithm have been demonstrated by computational experiments on two sets of widely used test instances. For the 46 challenging circle packing contest instances and the 24 widely-used NR instances, the algorithm can respectively improve 14 and 16 previous best-known records in a reasonable time.

There are two main directions for future research. On the one hand, the presented algorithm can be further improved by incorporating other advanced strategies. Possible improvements include the following: First, reduce the solution space by first ignoring several smaller disks and only looking for optimal packing pattern of the remaining larger disks. The smaller disks can be inserted into the holes after the larger disks have been placed into the container. This strategy was proposed in \cite{addis2008u} and had proved to be very useful. Second, test other Stochastic Local Search methods, such as Simulated Annealing used in \cite{muller2009}, Variable Neighborhood Search \citep{slsbook} and so on. Third, the proposed algorithm is a single-solution based method. It is possible to strengthen the robustness of the algorithm by employing some  population-based methods, like the Population Basin Hopping method proposed in \cite{grosso2007}.

On the other hand, the ideas behind the proposed algorithm can also be applied to other hard global optimization problems. Many real-world global optimization problems, such as the cluster optimization problem in computational chemistry \citep{blj2011} and the protein folding problem in computational biology \citep{huanglv2005},  have the same characteristics  as the unequal circle packing problem, i.e., they have both a continuous and combinatorial nature.  For these kinds of problems, it is possible to build a neighborhood structure on the set of local minima via appropriate perturbation moves, and then  to employ some advanced combinatorial optimization methods to systematically search for good local minima.

\section*{Acknowledgement}
We thank the anonymous reviewers whose detailed and valuable suggestions have significantly improved the manuscript.   We thank Eckard Specht for processing our data and publishing it on the Packomania website. This work was supported by National Natural Science Foundation of China (Grant No. 61100144, 61262011).



\begin{thebibliography}{00}
\bibitem[Addis et al.(2008a)]{addis2008e}
Addis, B., Locatelli, M., \& Schoen, F. (2008a). Disk packing in a square: a new global optimization approach. Informs Journal on Computing, 20, 516-524.

\bibitem[Addis et al.(2008b)]{addis2008u}
Addis, B., Locatelli, M., \& Schoen, F. (2008b). Efficiently packing unequal disks in a circle. Operations Research Letters, 36, 37-42.

\bibitem[Akeb et al.(2009)]{akeb2009}
Akeb, H., Hifi, M., \& M'Hallah, R. (2009). A beam search algorithm for the circular packing problem. Computers \& Operations Research, 36, 1513-1528.

\bibitem[Akeb and Hifi(2010)]{akeb2010a}
Akeb, H., \& Hifi, M. (2010). A hybrid beam search looking-ahead algorithm for the circular packing problem. Journal of Combinatorial Optimization, 20, 101-130.

\bibitem[Akeb et al. (2010)]{akeb2010b}
Akeb, H., Hifi, M., \& M'Hallah, R. (2010). Adaptive beam search lookahead algorithms for the circular packing problem. International Transactions in Operational Research, 17, 553-575.

\bibitem[Al-Mudahka et al. (2010)]{aimudahka2010}
Al-Mudahka, I., Hifi, M., \& M'Hallah, R. (2010). Packing circles in the smallest circle: an adaptive hybrid algorithm. Journal of the Operational Research Society, 62, 1917-1930.

\bibitem[Andreani et al. (2007)]{andreani2007}
Andreani, R., Birgin, E. G., Martínez, J. M., \& Schuverdt, M. L. (2007). On Augmented Lagrangian methods with general lower-level constraints. SIAM Journal on Optimization 18, 1286-1309.

\bibitem[Birgin et al.(2005)]{birgin2005}
Birgin, E. G., Mart\'{i}nez, J. M., \& Ronconi, D. P. (2005). Optimizing the packing of cylinders into a rectangular container: A nonlinear approach. European Journal of Operational Research, 160, 19-33.

\bibitem[Birgin and Sobral (2008)]{birgin2008}
Birgin, E. G., \& Sobral, F. N. C. (2008). Minimizing the object dimensions in circle and sphere packing problems. Computers \& Operations Research, 35, 2357-2375.

\bibitem[Birgin and Mart\'{i}nez (2009)]{birgin2009}
Birgin, E. G., \& Mart\'{i}nez, J. M. (2009). Practical Augmented Lagrangian Methods. In  C. A. Floudas \& P. M. Pardalos (Eds), Encyclopedia of Optimization (2nd ed., pp. 3013-3023), US: Springer.

\bibitem[Birgin et al.(2010)]{birgin2010}
Birgin, E. G., \& Gentil, J. M. (2010). New and improved results for packing identical unitary radius circles within triangles, rectangles and strips. Computers \& Operations Research, 37, 1318-1327.

\bibitem[Castillo et al.(2008)]{castillo2008}
Castillo, I., Kampas, F. J., \& Pint\'{e}r, J. D. (2008). Solving circle packing problems by global optimization: Numerical results and industrial applications. European Journal of Operational Research, 191, 786-80.

\bibitem[Fu et al.(2013)]{fu2013}
Fu, Z., Huang, W., \& L\"{u}, Z. (2013). Iterated tabu search for the circular open dimension problem. European Journal of Operational Research, 225, 236-243.

\bibitem[Glover and Laguna(1998)]{tabusearch}
Glover, F.,  \& Laguna, M. (1997). Tabu search.  Boston: Kluwer Academic Publishers.

\bibitem[Graham et al.(1998)]{graham1998}
Graham, R. L., Lubachevsky, B. D., Nurmela, K. J., \& \"{O}sterg{\aa}rd, P. R. J. (1998). Dense packings of congruent circles in a circle. Discrete Mathematics, 181, 139-154.

\bibitem[Grosso et al.(2007)]{grosso2007}
Grosso, A., Locatelli, M., \& Schoen, F. (2007). A population-based approach for hard global optimization problems based on dissimilarity measures. Mathematical Programming, 110, 373-404.

\bibitem[Grosso et al.(2010)]{grosso2010}
Grosso, A., Jamali, A. R., Locatelli, M., \& Schoen, F. (2010). Solving the problem of packing equal and unequal circles in a circular container. Journal of Global Optimization, 47, 63-81.

\bibitem[Hifi and M'Hallah(2004) ]{hifi2004}
Hifi, M., \& M'Hallah, R. (2004). Approximate algorithms for constrained circular cutting problems. Computers \& Operations Research, 31, 675-694.

\bibitem[Hifi and M'Hallah(2006)]{hifi2006}
Hifi, M., \& M'Hallah, R. (2006). Strip generation algorithms for constrained two-dimensional two-staged cutting problems. European Journal of Operational Research, 172, 515-527.

\bibitem[Hifi and M'Hallah(2007)]{hifi2007}
Hifi, M., \& M'Hallah, R. (2007). A dynamic adaptive local search algorithm for the circular packing problem. European Journal of Operational Research, 183, 1280-1294.

\bibitem[Hifi and M'Hallah(2008)]{hifi2008}
Hifi, M., \& M'Hallah, R. (2008). Adaptive and restarting techniques-based algorithms for circular packing problems. Computational Optimization and Applications, 39, 17-35.

\bibitem[Hifi and M'Hallah(2009)]{hifi2009}
Hifi, M., \& M'Hallah, R. (2009). A literature review on circle and sphere packing problems: Models and methodologies. Advances in Operations Research Volume 2009 (2009), Article ID 150624, doi:10.1155/2009/150624.

\bibitem[Hoos and St\"utzle(2005) ]{slsbook}
Hoos, H. H., \& St\"utzle, T. (2005). Stochastic local search: Foundations and applications. Morgan Kaufmann.

\bibitem[Huang et al.(2005)]{huang2005}
Huang, W., Li, Y., Akeb, H., \& Li, C. (2005). Greedy algorithms for packing unequal circles into a rectangular container. Journal of the Operational Research Society, 539-548.

\bibitem[Huang et al.(2006)]{huang2006}
Huang, W., Li, Y., Li, C., \& Xu, R. (2006). New heuristics for packing unequal circles into a circular container. Computers \& Operations Research, 33, 2125-2142.

\bibitem[Huang et al.(2006)]{huanglv2005}
Huang, W., Chen, M., \& L\"{u}, Z. (2006). Energy optimization for off-lattice protein folding. Physical Review E, 74, 41907.

\bibitem[Huang and Ye(2010)]{huang2010}
Huang, W., \& Ye, T. (2010). Greedy vacancy search algorithm for packing equal circles in a square. Operations Research Letters, 38, 378-382.

\bibitem[Huang and Ye(2011)]{huang2011}
Huang, W., \& Ye, T. (2011). Global optimization method for finding dense packings of equal circles in a circle. European Journal of Operational Research, 210, 474-481.

\bibitem[Huang et al. (2012a)]{huang2012a}
Huang, W, Fu, Z, Xu, R. (in press). Tabu search algorithm combined with global perturbation for packing arbitrary sized circles into a circular container,
Science China (Information Sciences), doi: http://dx.doi.org/10.1007/s11432-011-4424-3.

\bibitem[Huang et al. (2012b)]{huang2012b}
Huang, W. Zeng, Z., Xu, R., Fu, Z. (2012). Using iterated local search for efficiently packing unequal disks in a larger circle. Advanced Materials Research, 1477, 430-432.

\bibitem[L\"{u} and Huang(2008)]{lv2008}
L\"{u}, Z., \& Huang, W. (2008). PERM for solving circle packing problem. Computers \& Operations Research, 35, 1742-1755.

\bibitem[Liu and Nocedal(1989)]{lbfgs}
Liu, D. C., \& Nocedal, J. (1989). On the limited memory BFGS method for large scale optimization. Mathematical Programming, 45, 503-528.

\bibitem[L\'{o}pez and Beasley(2011)]{beasley2011}
L\'{o}pez, C. O., \& Beasley, J. E. (2011). A heuristic for the circle packing problem with a variety of containers. European Journal of Operational Research, 214, 512-525.

\bibitem[L\'{o}pez and Beasley(2012)]{beasley2012}
L\'{o}pez, C. O., \& Beasley, J. E. (2012). Packing unequal circles using formulation space search. Computers \& Operations Research. (accepted manuscript).

\bibitem[Louren\c{c}o et al.(2003)]{ils}
Louren\c{c}o, H., Martin, O.,\& St\"{u}tzle, T. (2003). Iterated local search. In F. Glover, \& G. Kochenberger (Eds.), Handbook of metaheuristics  (pp. 320-353). New York: Springer.

\bibitem[M\"{u}ller et al.(2009)]{muller2009}
M\"{u}ller, A., Schneider, J. J., \& Sch\"{o}mer, E. (2009). Packing a multidisperse system of hard disks in a circular environment. Physical Review E, 79, 21102.

\bibitem[Nurmela and \"{O}sterg{\aa}rd(1997)]{nurmela1997}
Nurmela, K. J., \&  \"{O}sterg{\aa}rd , P. R. J. (1997). Packing up to 50 equal circles in a square. Discrete \& Computational Geometry, 18, 111-120.

\bibitem[Specht(2013)]{specht2012}
Specht, E. (2013). Packomania website,  www.packomania.com.

\bibitem[Szab\'{o} et al.(2007)]{szabo2007}
Szab\'{o}, P. G., Mark\'{o}t, M. Cs., Csendes, T.,  Specht, E., Casado, L. G., \& Garc\'{i}a, I. (2007). New approaches to circle packing in a square. Springer.

\bibitem[Wang et al.(2002)]{wang2002}
Wang, H., Huang, W., Zhang, Q., \& Xu, D. (2002). An improved algorithm for the packing of unequal circles within a larger containing circle. European Journal of Operational Research, 141, 440-453.

\bibitem[Ye et al.(2011)]{blj2011}
Ye, T., Xu, R., \& Huang, W. (2011). Global optimization of binary lennard-jones clusters using three perturbation operators. Journal of Chemical Information and Modeling, 51, 572-577.

\bibitem[Zhang and Deng(2005)]{zhang2005}
Zhang, D., \& Deng, A. (2005). An effective hybrid algorithm for the problem of packing circles into a larger containing circle. Computers \& Operations Research, 32, 1941-1951.

\end{thebibliography}
\end{document}